\documentclass{amsart}
\usepackage{amsmath, amscd, amssymb, amsthm}
\usepackage{latexsym}
\usepackage{amsfonts}
\usepackage{graphicx}
\usepackage[all,cmtip]{xy}
\usepackage[colorlinks,linkcolor=blue,breaklinks=true,urlcolor=blue,citecolor=blue,anchorcolor=blue,pagebackref]{hyperref}%
\usepackage{geometry}
\newtheorem{theorem}{Theorem}
\newtheorem{lemma}{Lemma}

\newtheorem{remark}{Remark}

\newtheorem{conjecture}{Conjecture}

\renewcommand*\backref[1]{}
\renewcommand*\backrefalt[4]{ \ifcase #1 \or (cited on page #2) \else (cited on pages #2) \fi}

\newcommand{\be}{\begin{equation}}
\newcommand{\ee}{\end{equation}}
\newcommand{\bea}{\begin{eqnarray}}
\newcommand{\eea}{\end{eqnarray}}

\newcommand{\vs}{\vspace{0.5cm}}

\def\XXint#1#2#3{{\setbox0=\hbox{$#1{#2#3}{\int}$ }
\vcenter{\hbox{$#2#3$ }}\kern-.6\wd0}}

\begin{document}

\title[Lie hermitian manifolds with constant mixed curvature]{LIE  HERMITIAN MANIFOLDS WITH CONSTANT MIXED CURVATURE}

\author{Yuqin Guo}
\address{Yuqin Guo. School of Mathematical Sciences, Chongqing Normal University, Chongqing 401331, China}
\email{{1942747285@qq.com}}\thanks{Zhang is the corresponding author. Guo is supported by Chongqing graduate student research grant No. CYB25248.}

\author{Dongmei Zhang}
\address{Dongmei Zhang. School of Mathematical Sciences, Chongqing Normal University, Chongqing 401331, China}
\email{{2250825921@qq.com}}

\subjclass[2010]{53C55}
\keywords{complex manifold; Hermitian manifold; Bismut connection; mixed curvature}

\begin{abstract}
We study \text{Tang's conjecture} on compact Hermitian manifolds. We prove that this conjecture holds for two special classes of unimodular Lie groups: their associated Lie algebras are either almost abelian, or contain a J-invariant abelian ideal of codimension two. \end{abstract}

\maketitle


\markleft{Yuqin Guo and Dongmei Zhang}
\markright{Lie--Hermitian manifolds with constant mixed curvature}
\section{Introduction and statement of result}
Complex geometry harbours a long--standing open problem concerning classification of compact Hermitian manifolds endowed with constant holomorphic sectional curvature.
In the K\"ahler case, it is a well established fact that every compact (or complete) K\"ahler manifold of constant holomorphic sectional curvature belongs to the family of complex space forms.
In other words, their universal covers are one of the standard models: complex projective space$\,\mathbb{CP}^n$, complex Euclidean space $\mathbb{C}^n$, or complex hyperbolic space $\mathbb{CH}^n$.

For non-K\"ahler metrics $g$, the curvature tensor $R$ fails to satisfy the symmetry identities valid in the K\"ahler framework. Consequently, the holomorphic sectional curvature $H$ alone cannot recover full information encoded within $R$.
Nevertheless, the following statement remains widely conjectured to hold for compact manifolds.
\begin{conjecture}[Constant Holomorphic Sectional Curvature Conjecture] \label{conj1}
Let $(M^n,g)$ be a compact Hermitian manifold with $H$ equal to a constant $c$. If $c\neq 0$, then $g$ is K\"ahler and if $c=0$, then $g$ is Chern flat (namely $R=0$).
\end{conjecture}

Note that the compactness assumption is necessary for the validity of the conjecture, without it there are counterexamples in the noncompact case. Compact Chern flat manifolds have been classified by Boothby \cite{Boothby} as all the compact quotients of complex Lie groups equipped with left invariant metrics. Moreover, compact Levi-Civita flat Hermitian threefolds have been classified by Khan--Yang--Zheng in \cite{KYZ}, but for $n\geq4$, the classification remains an open question.

In complex dimension $n=2$, Conjecture \ref{conj1} was established through the work of Balas \cite{Balas}, Balas--Gauduchon \cite{BG}, Sato--Sekigawa \cite{SS}, and Apostolov--Davidov--Muskarov \cite{ADM}. In higher dimensions, an early result was obtained by Davidov--Grantcharov--Muskarov \cite{DGM}, who studied the conjecture for the natural Hermitian structures on twistor spaces. Chen--Chen--Nie \cite{CCN} verified Conjecture \ref{conj1} for compact locally conformally K\"ahler manifolds when $c\leq 0$. Huang--Wan \cite{HuangWan2026} subsequently removed the sign restriction on $c$.
Tang \cite{Tang} verified the conjecture for Chern-K\"ahler-like metrics, whose Chern curvature tensors satisfy all K\"ahler symmetries.
Further confirmations were obtained by Li--Zheng \cite{LZ} and Rao--Zheng \cite{RZ} for complex nilmanifolds and Bismut-K\"ahler-like metrics, respectively. Here, the Bismut-K\"ahler-like condition means that the Bismut curvature tensor satisfies all K\"ahler symmetries.

For real bisectional curvature, introduced by Yang--Zheng \cite{YangZ}, Zhou--Zheng \cite{ZhouZ} proved that compact Hermitian threefolds with vanishing real bisectional curvature are Chern flat.
Recently, Tang \cite{Tang2026RBC} resolved the Yang--Zheng conjecture by proving that, on a compact Hermitian manifold of complex dimension $n\geq 2$, constant real bisectional curvature must vanish and the metric is Chern flat.

Besides holomorphic sectional curvature, Chu--Lee--Tam introduced a new curvature invariant called mixed curvature, which generalizes the former notion.
For the precise definition, we follow Chu--Lee--Tam \cite{chu-lee-tam-1999}, who defined mixed curvature $\mathcal{C}_{\alpha,\beta}$ as a convex combination of the (first) Ricci curvature and holomorphic sectional curvature of the Chern connection.
\[
\mathcal{C}_{\alpha,\beta}(X) = \alpha \operatorname{Ric}(X, \overline{X}) + \beta H(X),
\]
where $\alpha, \beta \in \mathbb{R}$, $X \in T^{1,0}M$ is a type $(1,0)$ tangent vector of unit length, and $\operatorname{Ric}$ is the first Ricci curvature of the Chern connection, defined by $\operatorname{Ric}(X, \overline{X}) = \sum_{i=1}^n R(X, \overline{X}, e_i, \overline{e}_i)$ where $e$ is any local unitary frame.
Throughout this paper, we will always assume that $(\alpha, \beta) \neq (0,0)$.

In the above notation, we have $\mathcal{C}_{1,0} = \operatorname{Ric}$ and $\mathcal{C}_{0,1} = H$. When the metric $g$ is K\"ahler, the mixed curvature encompasses several interesting curvature conditions.
Specifically, $\mathcal{C}_{1,1}$ corresponds to the $\operatorname{Ric}^+(X, \overline{X})$ introduced by Ni \cite{ni-2004}, and $\mathcal{C}_{1,-1}$ corresponds to the orthogonal Ricci curvature $\operatorname{Ric}^\perp(X, \overline{X})$ introduced by Ni--Zheng \cite{ni-zheng-2008}.
Additionally, if $k$ is an integer between $1$ and $n$, then $\mathcal{C}_{k-1,n-k}$ is closely related to the $k$-Ricci curvature $Ric_k$ introduced by Ni \cite{ni-2005}.

For compact K\"ahler manifolds with $\mathcal{C}_{\alpha,\beta} \geq 0$ or $\mathcal{C}_{\alpha,\beta} \leq 0$, there have been important progress made in recent years.
For example, Yang \cite{yang-2007} proved that a compact K\"ahler manifold with $\mathcal{C}_{0,1} > 0$ must be projective and rationally connected, confirming a conjecture by Yau \cite{yau-1978}.
Ni \cite{ni-2004} generalized this result, showing that the same holds if $Ric_k > 0$ for some $1 \leq k \leq n$.
Matsumura \cite{matsumura-1983}, \cite{matsumura-1985} established structure theorems for projective manifolds with $\mathcal{C}_{0,1} \geq 0$, while Chu--Lee--Tam \cite{chu-lee-tam-1999} and others proved that compact K\"ahler manifolds with $\mathcal{C}_{\alpha,\beta} > 0$, where $\alpha > 0$ and $3\alpha + 2\beta \geq 0$, must be projective and simply-connected.
Zhang--Zhang \cite{zhang-zhang-2010} extended Yang's result to the quasi-positive case, confirming a conjecture by Yang.
Tang \cite{tang-2012} demonstrated the projectivity of compact K\"ahler manifolds with quasi-positive $\mathcal{C}_{\alpha,\beta}$.
More recently, Chu--Lee--Zhu \cite{chu-lee-zhu-2018} established a structure theorem for compact K\"ahler manifolds with $\mathcal{C}_{\alpha,\beta} \geq 0$.
In the non-positive case, Wu--Yau \cite{wu-yau-2014} confirmed a conjecture of Yau that projective K\"ahler manifolds with $\mathcal{C}_{0,1} < 0$ must have ample canonical line bundles.
Tosatti--Yang \cite{tosatti-yang-2013} showed that the projectivity assumption can be removed in Wu-Yau's Theorem, and Chu--Lee--Tam \cite{chu-lee-tam-1999} proved that compact K\"ahler manifolds with $\mathcal{C}_{\alpha,\beta} < 0$ and $\alpha \geq 0$, $\beta \geq 0$, must have ample canonical bundles.

Most existing results for mixed-curvature invariants focus predominantly on the K\"ahler setting; comparatively little is understood for non-K\"ahler geometries. In \cite{tang-2012}, Tang formulated the following conjecture:

\begin{conjecture}[Tang's conjecture] \label{conj2}
Let $(M^n, g)$ be a compact Hermitian manifold with $\mathcal{C}_{\alpha,\beta}=c$, where $c$ is a constant. If $c \neq 0$, then $g$ is K\"ahler.
\end{conjecture}

Note that in the conjecture one may assume that $\beta \neq 0$. The $\beta = 0$ case corresponds to (first) Chern Ricci curvature, where the metric is automatically K\"ahler when $\mathcal{C}_{1,0} = c \neq 0$, as in this case the K\"ahler form is a constant multiple of the Chern Ricci form, which is always closed.

For the $c = 0$ case, the vanishing of the mixed curvature $\mathcal{C}_{\alpha,\beta}$ in general does not imply $R = 0$, Tang illustrated this phenomenon in \cite{tang-2012} via examples of isosceles Hopf manifolds.

In complex dimension $n = 2$, Tang \cite{tang-2012} established that every compact Hermitian surface of constant mixed curvature must either be K\"ahler or an isosceles Hopf surface.
For dimensions $n \geq 3$, Tang also obtained results covering locally conformally K\"ahler manifolds and Chern-K\"ahler-like manifolds.
Chen--Zheng \cite{chen-zheng-2025} verified Conjecture \ref{conj2} across multiple families of special Hermitian manifolds: complex nilmanifolds, solvmanifolds with complex commutators, almost-Abelian Lie groups, and Lie algebras containing a J-invariant Abelian ideal of codimension two.
They also verified Conjecture \ref{conj2} for all compact balanced three-dimensional manifolds whose Bismut connection has parallel torsion.

One major obstacle to resolving Conjecture \ref{conj1} lies in algebraic complexity of curvature tensors on Hermitian manifolds.
Away from the K\"ahler domain, curvature tensors for Chern, Levi-Civita and Bismut connections lose the familiar K\"ahlerK symmetry properties.
As a consequence, holomorphic sectional curvature only governs the symmetrized part of curvature tensor and cannot fully determine the complete curvature data.

This motivates focusing on the subclass of Lie-Hermitian manifolds---namely, compact Hermitian spaces of the form $M = G/\Gamma$, where $G$ is a Lie group and $\Gamma \subset G$ is a discrete cocompact subgroup, with both the complex structure $J$ and the metric $g$ lifted to $G$ as left-invariant objects as a useful test case.
There exists a wealth of such triples $(G, J, g)$ for which $J$ is not bi-invariant, meaning that  $G$ need not be a complex Lie group.
Within this framework, the manifold is termed a complex nilmanifold when $G$ is nilpotent, and a complex solvmanifold when $G$ is solvable.

It follows readily from the definition of mixed curvature that it is a combination of the holomorphic sectional curvature and the Ricci curvature.
Hence Conjecture \ref{conj2} generalizes Conjecture \ref{conj1}. It is natural to specialize Conjecture \ref{conj2} to this particular Lie-Hermitian setting and investigate the corresponding weaker version of the problem.
Observe that the existence of a compact quotient forces $G$ to be unimodular. Moreover, according to a classical theorem of Hano \cite{hano-1957}, every unimodular K\"ahler Lie algebra must be flat.
Consequently, when Conjecture \ref{conj1} is restricted to the Lie-Hermitian realm, the problem reduces to the following\ldots

The main purpose of this short article is to explore Conjecture \ref{conj2} for two special types of Lie-Hermitian manifolds: when the Lie algebra ${\mathfrak g}$ of $G$ is almost abelian (namely, it contains an abelian idael of codimension $1$), or when ${\mathfrak g}$ contains an abelian ideal ${\mathfrak a}$ of codimension $2$ satisfying $J{\mathfrak a}={\mathfrak a}$, which is a natural generalization to the almost abelian case.
In the almost abelian case, we are able to give a complete answer to Conjecture \ref{conj2}.
\begin{theorem} \label{thm1}
Let $(M^n, g)$ be a Lie-Hermitian manifold with universal cover $(G, J, g)$. Assume that the Lie algebra $\mathfrak{g}$ of $G$  is unimodular and almost abelian.
If $\mathcal{C}^{r}_{\alpha,\beta} = c$, then $c = 0$. If $\mathcal{C}^{r}_{\alpha,\beta} = 0$ and $\beta \neq 0$, then $g$ must be K\"ahler flat.
\end{theorem}
Besides the results for the Riemannian connection, we also present counterparts for the Bismut connection.
\begin{theorem} \label{thm2}
Let $(M^n, g)$ be a Lie-Hermitian manifold with universal cover $(G, J, g)$. Assume that the Lie algebra $\mathfrak{g}$ of $G$  is unimodular and almost abelian.
If $\mathcal{C}^{b}_{\alpha,\beta} = c$, then $c = 0$. If $\mathcal{C}^{b}_{\alpha,\beta} = 0$ and $\beta \neq 0$, then $g$ must be K\"ahler flat.
\end{theorem}

In the ${\mathfrak g}$ contains an abelian ideal ${\mathfrak a}$ of codimension $2$ satisfying $J{\mathfrak a}={\mathfrak a}$ case, we are able to give a complete answer to Conjecture \ref{conj2}.
\begin{theorem} \label{thm3}
Let $(M^n, g)$ be a Lie-Hermitian manifold with universal cover $(G, J, g)$. Assume that the Lie algebra $\mathfrak{g}$ of $G$ is unimodular and contains a $J$-invariant abelian ideal of codimension $2$. If $\mathcal{C}^{r}_{\alpha,\beta} = c$, then $c = 0$. If $\mathcal{C}^{r}_{\alpha,\beta} = 0$ and $\beta \neq 0$, then $g$ must be Chern flat.
\end{theorem}

\vspace{0.3cm}
\section{Preliminaries }
In this section, let us set up the notations. Let $(M^n,g)$ be a Hermitian manifold and denote by $\omega$ the  K\"ahler form associated with $g$. We will follow the notations of \cite{ChenZ25}. Fix any point $p\in M$, let $\{ e_1, \ldots , e_n\}$ be a local frame of type $(1,0)$ complex tangent vector fields on $M^n$ around $p$, and let $\varphi=\, ^t\! (\varphi_1,...,\varphi_n)$ denote the dual coframe. Write $g=\langle, \rangle$ and extend it bi-linearly over ${\mathbb C}$.

Let $\nabla$, $\nabla^c$, and $\nabla^b$ denote the Levi-Civita (Riemannian), Chern, and Bismut connection, respectively. Denote by $R$ the curvature of $\nabla$, by $T^c=T$ and $R^c$ the torsion and curvature of $\nabla^c$, and by $T^b$ and $R^b$ the torsion and curvature of $\nabla^b$. Under any local unitary frame $e$, the components of $T^c$ are given by
$$T^c(e_i,\bar{e}_j)=0, \ \ \ T^c(e_i,e_j)=\sum^n_{k=1}T^k_{ij}e_k,$$
where $1\leq i,j\leq n$ and $T^k_{ij}=-T^k_{ji}$. For the Chern connection $\nabla^c$, let $\theta$ and $\Theta$ denote the matrix of connection and curvature, respectively, with respect to a local unitary frame $e$. Let $\tau$ be the column vector of the Chern torsion forms under $e$, namely, $\tau_j=\frac{1}{2} \sum_{i,k}T^j_{ik}\varphi_i\wedge \varphi_k$. Then the structure equations and Bianchi identities are
\begin{equation*}
\label{structure-c} \left\{ \begin{array}{llll} d \varphi = - \ ^t\!\theta \wedge \varphi + \tau,   \\
d  \theta = \theta \wedge \theta + \Theta.  \\
d \tau = - \ ^t\!\theta \wedge \tau + \ ^t\!\Theta \wedge \varphi, \\
d  \Theta = \theta \wedge \Theta - \Theta \wedge \theta. \end{array} \right.
\end{equation*}
Similarly, denote by $\theta^b$ and $\Theta^b$ the matrix of connection and curvature of $\nabla^b$ under $e$. Then
\begin{equation*}
\label{structure-s}
\Theta^b = d\theta^b -\theta^b \wedge \theta^b.
\end{equation*}
Let $\gamma = \nabla^b-\nabla^c$ be the tensor, and for simplicity we will also write $\gamma = \theta^b -\theta$ for its matrix representation under $e$. Then by \cite{YZ} we have
$$\gamma e_i = \sum_j \gamma_{ij} e_j = \sum_{j,k} \big( T^j_{ik}\varphi_k - \overline{T^i_{jk} } \, \overline{\varphi}_k \big) e_j.$$
Moreover, following the notations in \cite{YZ,YZ1}, the Levi-Civita connection is given by
\begin{equation*}
 \nabla e_i = \sum_j \{ \big( \theta_{ij} +\frac{1}{2}\gamma_{ij}\big)  e_j  +  \overline{\beta_{ij}} \,\overline{e}_j \} , \ \ \ \mbox{where} \ \ \  \beta_{ij} = \frac{1}{2} \sum_k \overline{T^k_{ij}} \varphi_k.
\end{equation*}

\vspace{0.1cm}
Denote by $R^D$ the curvature tensor of the canonical metric connection $D=D^r_s$. The holomorphic sectional curvature of $D=D^r_s$, denoted by $H^D$, is defined as
$$H^D(X)=R^D_{X\overline{X}X\overline{X}}/|X|^4$$
for any nonzero tangent vector X of type (1,0). By definition,  $H^D=c$ for a constant $c$ means $$R^D_{X\overline{X}X\overline{X}}=c|X|^4, \ \ \ \ $$
for all type $(1,0)$ tangent vector $X$. Under any local unitary frame $e$, the components of $R^D$ are denoted by $R^D_{i\bar{j}k\bar{\ell}}=R^D(e_i, \overline{e}_j , e_k, \overline{e}_{\ell})$, and we have
\begin{equation}\label{HS}
 \widehat{R}^D_{i\bar{j}k\bar{\ell}} = \frac{c}{2}(\delta_{ij}\delta_{k\ell } + \delta_{i\ell }\delta_{kj}).
\end{equation}
where
\begin{equation}\label{RD}
 \widehat{R}^D_{i\bar{j}k\bar{\ell}}  = \frac{1}{4} \big( R^D_{i\bar{j}k\bar{\ell}}+  R^D_{k\bar{j}i\bar{\ell}} +  R^D_{i\bar{\ell}k\bar{j}} +  R^D_{k\bar{\ell}i\bar{j}} \big)
\end{equation}
is the symmetrization of $R^D$. Here, $D_{0}^{1}=\nabla^{c}$, $\nabla^{b}=D_{0}^{-1}$. The Riemannian connection is a linear combination of the Chern connection and the Bismut connection, and the relation is given below.
$$ \nabla= \frac{1}{2} \nabla^c + \frac{1}{2}\nabla^b.$$

Next, let us recall some basic properties for {\em Lie-Hermitian manifolds,} which are compact Hermitian manifolds whose universal cover is $(G,J,g)$, where $G$ is a Lie group equipped with a left-invariant complex structure $J$ and a compatible left-invariant metric $g$. Let ${\mathfrak g}$ denote the Lie algebra of $G$. Next, let us recall some basic properties for Lie-Hermitian manifolds. Let $G$ be a connected, simply-connected, even-dimensional Lie group, and $\mathfrak{g}$ its Lie algebra. Left-invariant complex structures and compatible left-invariant metrics on $G$ correspond to almost complex structures $J$ and compatible inner products $g=\langle\,,\rangle$ on $\mathfrak{g}$, such that $J$ is integrable, namely,
\begin{equation}\label{integrability}
[x,y]-[Jx,Jy]+J[Jx,y]+J[x,Jy]=0,\qquad \forall\,x,y\in\mathfrak{g}.
\end{equation}
Extend $J$ and $\langle\,,\rangle$ linearly over $\mathbb{C}$ to the complexification $\mathfrak{g}^{\mathbb{C}}$.
A left-invariant complex structure on $G$ corresponds to a complex structure on ${\mathfrak g}$ (which we also denote by $J$ for convenience). A compatible  left-invariant metric on $G$ corresponds to an inner product $g=\langle, \rangle$ on ${\mathfrak g}$ satisfying $\langle Jx, Jy \rangle = \langle x, y \rangle$ for all $x, y \in {\mathfrak g}$. Here again for convenience we use the same letter to denote the corresponding metric on the Lie algebra. Let $ \mathfrak{g}^{\mathbb{C}} $ be the complexification of $ \mathfrak{g} $, and define its $(1,0)$-part by
\[
\mathfrak{g}^{1,0} := \{ x - \sqrt{-1}Jx \mid x \in \mathfrak{g} \} \subseteq \mathfrak{g}^{\mathbb{C}}.
\]
The integrability condition \eqref{integrability} ensures that $ \mathfrak{g}^{1,0} $ is a complex Lie subalgebra of $ \mathfrak{g}^{\mathbb{C}} $. We extend $J$ and $ g = \langle\, , \, \rangle $ complex bilinearly to $ \mathfrak{g}^{\mathbb{C}} $, and assume $ \dim_{\mathbb{R}} G = 2n $. A \emph{unitary frame} of $\mathfrak{g}$ is a basis $ e = \{ e_1, \ldots, e_n \} \subset \mathfrak{g}^{1,0} \cong \mathbb{C}^n $ such that $\langle e_i, \overline{e}_j \rangle = \delta_{ij}$ for any $1 \leq i,j \leq n.$
Let $ \varphi = \{ \varphi_1, \ldots, \varphi_n \} $ be the dual coframe of $e$, which means a basis of $ (\mathfrak{g}^{1,0})^* $ satisfying $ \varphi_i(e_j) = \delta_{ij} $. Following the conventions in \cite{VYZ, YZ1, ZZ1}, the structure constants are given by
\begin{equation*}
C^j_{ik} := \varphi_j([e_i, e_k]), \ \ \ \quad D^j_{ik} := \overline{\varphi}_i([\overline{e}_j, e_k]), \ \ \ \quad \forall\, 1 \leq i,j \leq n,
\end{equation*}
which is equivalent to
\begin{equation} \label{CandD}
[e_i, e_j] = \sum_k C^k_{ij} e_k, \ \ \ \ \ \ \ \ [e_i, \overline{e}_j] = \sum_k \big( \overline{D^i_{kj}} e_k - D^j_{ki} \overline{e}_k \big), \ \ \ \ \ \forall \ 1\leq i,j\leq n.
\end{equation}
In dual terms, these relations are also equivalent to the (first) structure equation:
\begin{equation}\label{eq:structure}
d\varphi_i = -\sum_{j,k=1}^n \big( \frac{1}{2} C^i_{jk}\varphi_j\wedge \varphi_k  + \overline{D^j_{ik}} \varphi_j \wedge \overline{\varphi}_k \big).
\end{equation}
The first Bianchi identity, which coincides with the Jacobi identity, takes the form
\begin{equation*}
\left\{ \begin{aligned}  \sum_{r=1}^n \big( C^r_{ij}C^{\ell}_{rk} + C^r_{jk}C^{\ell}_{ri} + C^r_{ki}C^{\ell}_{rj} \big) \ = \ 0,   \hspace{3.2cm}\\
 \sum_{r=1}^n \big( C^r_{ik}D^{\ell}_{jr} + D^r_{ji}D^{\ell}_{rk} - D^r_{jk}D^{\ell}_{ri} \big) \ = \ 0,  \hspace{3cm} \\
 \sum_{r=1}^n \big( C^r_{ik}\overline{D^{r}_{j\ell }} - C^j_{rk}\overline{D^{i}_{r\ell }} + C^j_{ri}\overline{D^{k}_{r\ell }}  - D^{\ell}_{ri}\overline{D^{k}_{jr }} + D^{\ell}_{rk} \overline{ D^{i}_{jr }}  \big) \ = \ 0 ,  \end{aligned} \right.
\end{equation*}
for any $1\leq i,j,k,\ell \leq n$. Under any unitary frame $e$, the Chern connection form $\theta_{ij}$, the Chern torsion components $T^j_{ik}$, and the Chern curvature components $R_{i\bar{j}k\bar{\ell}}$ take particularly simple forms (\cite{LZ}):
\begin{eqnarray}
 &&  \theta_{ij} \ = \ \sum_{k=1}^n \big( D^j_{ik}\varphi_k  -\overline{D^i_{jk}} \overline{\varphi_k} \big), \ \ \ \ \ \
  T^j_{ik}  \ = \    - C^j_{ik} -   D^{j}_{i k} + D^{j}_{k i} , \ \ \ \ \forall \ 1\leq i,j,k \leq n. \label{torsion} \\
 && R^c_{i\bar{j}k\bar{\ell}} \ = \   \sum_{s=1}^n \big( D^s_{ki}\overline{D^s_{\ell j}} - D^{\ell}_{si}\overline{D^k_{s j}} - D^j_{si}\overline{D^k_{ \ell s}} - \overline{D^i_{sj}} D^{\ell}_{k s} \big), \  \ \ \ \ \ \forall\ 1\leq i,j,k,\ell \leq n. \label{curvature}
\end{eqnarray}
Under any unitary basis e of ${\mathfrak g}$, the components of mixed  curvature tensor $R^{b}$ of a Lie Hermitian manifold are given by
\begin{lemma}[\cite{chen-zheng-2025}]\label{lemma1}
Let $(M^n, g)$ be a Hermitian manifold. Then $\mathcal{C}_{\alpha,\beta} = c \Longleftrightarrow$
\[
\quad 4\beta\widehat{R}_{i\bar{j}k\bar{\ell}} + \alpha\bigl(R_{i\bar{j}}g_{k\bar{\ell}} + R_{k\bar{j}}g_{i\bar{\ell}} + R_{i\bar{\ell}}g_{k\bar{j}} + R_{k\bar{\ell}}g_{i\bar{j}}\bigr) = 2c\bigl(g_{i\bar{j}}g_{k\bar{\ell}} + g_{i\bar{\ell}}g_{k\bar{j}}\bigr), \quad \forall\ 1 \leq i,j,k,\ell \leq n,
\]
where $R_{i\bar{j}} = \operatorname{Ric}(e_i, \bar{e}_j)$ denotes the components of the first Chern Ricci curvature. In particular, when the frame $e$ is unitary, we have $R_{i\bar{j}} = \sum_{s=1}^n R_{i\bar{j}s\bar{s}}$ and
\begin{align}
\label{eq:lem-eq1}
\quad &\beta R_{i\bar{i}i\bar{i}} + \alpha R_{i\bar{i}} = c, \quad \forall\ 1 \leq i \leq n, \\
\label{eq:lem-eq2}
\quad &4\beta\widehat{R}_{i\bar{i}k\bar{k}} + \alpha\bigl(R_{i\bar{i}} + R_{k\bar{k}}\bigr) = 2c, \quad \forall\ 1 \leq i \neq k \leq n.
\end{align}
\end{lemma}

\vspace{0.5cm}

\section{The almost abelian case}
In this section, we study Lie-Hermitian manifolds whose universal covers are modeled on \emph{almost abelian Lie algebras}, i.e., real Lie algebras \( \mathfrak{g} \) admitting an abelian ideal \( \mathfrak{a} \subset \mathfrak{g} \) of codimension $1$. Throughout this section, we adopt the notations of \cite{GuoZ} and assume that ${\mathfrak g}$ is almost abelian.

Let $(\mathfrak{g}, J, g)$ be a Hermitian Lie algebra that is almost abelian, with ${\mathfrak a}\subset {\mathfrak g}$ an abelian ideal of codimension $1$. Then
${\mathfrak a}_J:={\mathfrak a}\cap J {\mathfrak a}$  is a $J$-invariant ideal of codimension $2$ in ${\mathfrak g}$. A unitary basis  $e$ of ${\mathfrak g}^{1,0}$ is called an {\em admissible frame} for $({\mathfrak g}, J, g)$  if ${\mathfrak a}_J$ is spanned
by the real and imaginary parts of $e_i$ for $2\leq i\leq n$, while
$$ {\mathfrak a} = \mbox{span}_{\mathbb R}\{ \sqrt{-1}(e_1-\overline{e}_1); \ (e_i+\overline{e}_i), \, \sqrt{-1}(e_i-\overline{e}_i); \ 2\leq i\leq n \} .$$
From now on, we will assume that $e$ is an admissible frame. Following \cite{GuoZ}, since ${\mathfrak a} $ is abelian, we have
$$ [e_i, e_j]= [ e_i, \overline{e}_j] = [e_1-\overline{e}_1, e_i]= 0, \ \ \ \ \forall \ 2\leq i,j\leq n.$$
Using this and equation (\ref{CandD}), we conclude that
\begin{equation*}
C^{\ast}_{ij} = D^j_{\ast i} = D^1_{\ast i} = C^{\ast}_{1i}+\overline{D^i_{\ast 1}} =0, \ \ \ \ \forall \ 2\leq i,j\leq n.
\end{equation*}
Moreover, since ${\mathfrak a}$ is an ideal, $[e_1+\overline{e}_1, {\mathfrak a} ] \subseteq {\mathfrak a}$, which implies that $C^1_{1i}=0$ and $D^1_{11} = \overline{D^1_{11}}$. Combining these results, we know that the only possibly non-trivial components of $C$ and $D$ are
\begin{equation} \label{CandD-aala}
D^1_{11}=\lambda \in {\mathbb R}, \ \ \ D^1_{i1}=v_i \in {\mathbb C}, \ \ \ D^j_{i1}=A_{ij}, \ \ \ C^{j}_{1i} = - \overline{A_{ji}}, \ \ \ \ \ \  2\leq i,j\leq n.
\end{equation}
In other words, for any almost abelian Hermitian Lie algebra and under any admissible frame, the structure constants are determined by a real number $\lambda \in {\mathbb R}$, a complex column vector $v\in {\mathbb C}^{n-1}$, and a complex $(n-1)\times (n-1)$ matrix $A=(A_{ij})_{2\leq i,j\leq n}$. With respect to the dual coframe
$\varphi$, the structure equation (\ref{eq:structure}) takes the form
$$
\begin{cases}
d\varphi_1 = - \lambda \,\varphi_1\wedge \overline{\varphi}_1 ,\vspace{0.2cm} \\
d\varphi_i = - \overline{v}_i \, \varphi_1\wedge \overline{\varphi}_1 +  \sum^n_{j=2} \overline{A_{ij}} \,(\varphi_1 + \overline{\varphi}_1)\wedge \varphi_j , \ \ \ 2\leq i\leq n.
\end{cases}
$$
To prove Theorem \ref{thm1} for almost abelian Lie algebras, we will need the following characterization result from \cite[Proposition 7 and Lemma 6]{GuoZ} for almost abelian Hermitian Lie algebras.

\begin{lemma}[\cite{YZ}]\label{lemma2}
Let $(M^n,g)$ be a Hermitian manifold and $e$ be a local unitary frame.
The Chern and Levi--Civita curvature components under $e$ are related by
\begin{equation}\label{eq:4}
\begin{cases}
T_{ik,\bar{j}}^{\ell}=R_{k\bar{j}i\bar{\ell}}^{c}-R_{i\bar{j}k\bar{\ell}}^{c},\\[4pt]
R_{i\bar{j}k\bar{\ell}}=R_{i\bar{j}k\bar{\ell}}^{c}+\dfrac12\big(T_{ik,\bar{j}}^{\ell}+\overline{T_{j\ell,\bar{i}}^{k}}\big)
+\dfrac14\big(T_{ik}^{s}\overline{T_{j\ell}^{s}}-T_{is}^{\ell}\overline{T_{js}^{k}}-T_{ks}^{j}\overline{T_{\ell s}^{i}}\big),\\[4pt]
\end{cases}
\end{equation}
\end{lemma}

\begin{lemma}[\cite{LZ25}]
Let $(M^n,g)$ be a Hermitian manifold and $e$ be a local unitary frame.
The symmetrization of Levi--Civita curvature has components
\[
\widehat{R}_{i\bar{j}k\bar{\ell}}
=\widehat{R}^{c}_{i\bar{j}k\bar{\ell}}
-\frac18\sum_{s=1}^{n}
\big(T_{is}^{j}\overline{T_{\ell s}^{k}}
+T_{ks}^{\ell}\overline{T_{js}^{i}}
+T_{is}^{\ell}\overline{T_{js}^{k}}
+T_{ks}^{j}\overline{T_{\ell s}^{i}}\big),
\quad \forall\,1\le i,j,k,\ell\le n.
\]

In particular, we have
\begin{equation}\label{eq:5}
\begin{cases}
R_{i\bar{i}i\bar{i}}
=R^{c}_{i\bar{i}i\bar{i}}
-\dfrac12\displaystyle\sum_{s=1}^{n}\big|T_{is}^{i}\big|^{2},
\quad \forall\,1\le i\le n.\\[8pt]
\widehat{R}_{i\bar{i}k\bar{k}}
=\widehat{R}^{c}_{i\bar{i}k\bar{k}}
-\dfrac18\displaystyle\sum_{s=1}^{n}
\big(\big|T_{ks}^{i}\big|^{2}+\big|T_{is}^{k}\big|^{2}
+2\operatorname{Re}\big\{T_{is}^{i}\,\overline{T_{ks}^{k}}\big\}\big),
\quad \forall\,1\le i,k\le n.
\end{cases}
\end{equation}
\end{lemma}

\begin{lemma} [\cite{GuoZ}]\label{lemma3}
Let ${\mathfrak g}$ be an almost abelian Lie algebra equipped with a Hermitian structure $(J, g)$, and let $e$ be an admissible frame of ${\mathfrak g}$. Then the following hold:
\begin{enumerate}
\item   ${\mathfrak g}$ is unimodular $\ \Longleftrightarrow \ $ $\lambda + \mbox{tr}(A) + \overline{\mbox{tr}(A)} = 0$;
\item $g$ is K\"ahler  $\ \Longleftrightarrow \ $ $v=0$ and $A+A^{\ast}=0$.
\item $g$  is Chern flat $\ \Longleftrightarrow \ $ $\lambda=0$, $v=0$, and $[A, A^{\ast}]=0$.
 \end{enumerate}
 Here $A^{\ast}$ stands for the conjugate transpose of the matrix $A$.
\end{lemma}

As stated in \cite{GuoZ}, it is known that for any admissible frame $e$, the only possibly non-zero components of the Chern torsion and curvature are
\begin{eqnarray}
\begin{cases}
T^1_{1i} = v_i, \ \ \ \ T^j_{1i} = A_{ij} + \bar{A}_{ji}, \\
R^{c}_{1\bar{1}1\bar{1}} = -2\lambda^2 - |v|^2,  \ \ \quad (R^{c}_{1\bar{1}i\bar{1}}) = -A^*v, \\
(R^{c}_{1\bar{1}i\bar{j}}) = vv^* + [A, A^*] - \lambda(A + A^*),
\end{cases} \label{almostcur}
\end{eqnarray}
for any $2 \leq i,j \leq n$.

\vspace{0.5cm}
\begin{proof}[{\bf Proof of Theorem \ref{thm1}}]
Let $(\mathfrak{g}, J, g)$ be a unimodular almost abelian Lie algebra with a Hermitian structure. Assume that the Riemann connection of $g$ has constant mixed curvature: $\mathcal{C}_{\alpha,\beta} = c$, where $\beta \neq 0$. Let  $e$ be an admissible frame and write $S=A+A^*$. After a unitary change of $\{e_2,...,e_n\}$ if necessary, we may assume that $S$ is a real diagonal matrix. Since $\mathcal{C}_{\alpha,\beta} = c$, by (\ref{eq:lem-eq1}),(\ref{eq:lem-eq2}),We perform classification and obtain the following equation:
\begin{equation}\label{eq:system-mixed}
\begin{cases}
4\beta\widehat{R}_{i\bar{i}k\bar{k}} + \alpha\bigl(R_{i\bar{i}} + R_{k\bar{k}}\bigr) = 2c, \quad \forall \ 2 \leq i \neq k \leq n. \\
4\beta\widehat{R}_{1\bar{1}i\bar{i}} + \alpha\bigl(R_{1\bar{1}}+ R_{i\bar{i}}\bigr) = 2c, \ \quad \forall \ 2 \leq i \leq n.\\
\beta R_{i\bar{i}i\bar{i}} + \alpha R_{i\bar{i}} = c, \ \ \ \ \ \ \ \ \ \ \ \ \ \ \ \ \quad \forall \ 2 \leq i \leq n.\\
\beta R_{1\bar{1}1\bar{1}} + \alpha R_{1\bar{1}} = c.
\end{cases}
\end{equation}
From equations (\ref{eq:4}), (\ref{eq:5}) and (\ref{almostcur}), one can derive by computation that:
\begin{equation*}
\hspace{2em}
\begin{cases}
R_{1\bar{1}1\bar{1}} = -2\lambda^2 - \frac{3}{2}|v|^2, \quad {R}_{i\bar{i}i\bar{i}} = -\frac{1}{2}|S_{ii}|^2, \quad \widehat{R}_{i\bar{i}k\bar{k}} = -\dfrac{1}{4} S_{ii}S_{kk}, \quad \forall \ 2 \leq i \neq k \leq n. \\
\widehat{R}_{1\bar{1}i\bar{i}} = \dfrac{1}{4} \big(|v_{i}|^{2} + [A, A^*]_{ii} - \lambda S_{ii}\big) -\frac{1}{8}\big(|v_{i}|^{2}+|S_{ii}|^2\big), \quad \forall \ 2 \leq i,j \leq n. \\
R_{1\bar{1}} = -2\lambda^2 - \frac{3}{2}|v|^2, \quad {R}_{i\bar{i}} =\frac{1}{4}\sum_{m}\big(|S_{im}|^2-|S_{mi}|^2\big), \quad \forall \ 2 \leq i \leq n. \\
\end{cases}
\end{equation*}
Substituting the above formula into equation (\ref{eq:system-mixed}), we get
\begin{align}
c &= \alpha \left\{ -2\lambda^2 - \frac{3}{2}|v|^2 \right\}
+ \beta \left\{ -2\lambda^2 - \frac{3}{2}|v|^2 \right\}, \label{eq:12}\\[6pt]
c &= \alpha \left\{ \frac14\sum_{m}\big(|S_{im}|^2 - |S_{mi}|^2\big) \right\}
+ \beta \left\{ -\frac12|S_{ii}|^2 \right\},\quad 2\le i\le n. \label{eq:13}\\[6pt]
2c &= \alpha\left\{ \frac14\sum_{m}\big(|S_{im}|^2-|S_{mi}|^2\big)
-2\lambda^2-\frac23|v|^2 \right\} \notag\\
&\quad + \beta\left\{ [A,A^*]_{ii}-\lambda S_{ii}
-\frac12|S_{ii}|^2+\frac12|v_i|^2 \right\},\quad 2\le i\le n. \label{eq:14}\\[6pt]
2c &= \alpha\left\{ \frac14\sum_{m}\big(|S_{im}|^2-|S_{mi}|^2
+|S_{km}|^2-|S_{mk}|^2\big) \right\}
+ \beta\left\{-S_{ii}S_{kk} \right\},\quad 2\le i\neq k\le n. \label{eq:15}
\end{align}
Adding equations (\ref{eq:12}) and (\ref{eq:13}), we obtain
$$2c=\alpha\left\{-2\lambda^2 - \frac{3}{2}|v|^2+\frac14\sum_{m}\big(|S_{im}|^2 - |S_{mi}|^2\big)\right\}
+\beta\left\{-2\lambda^2 - \frac{3}{2}|v|^2-\frac12|S_{ii}|^2 \right\}.$$
By subtracting equation (\ref{eq:14}) from the above equation, we obtain the following key identity.
$$0=\beta\left\{-2\lambda^2 - \frac{3}{2}|v|^2-\frac12|v_{i}|^2 -[A,A^*]_{ii}+\lambda S_{ii}\right\}$$
Summing both sides of the above identity over $i$, and recalling that the Lie algebra $\mathfrak{g}$ is unimodular, i.e., $\lambda=-\operatorname{tr}(S)$, we substitute this relation and obtain the following result.
$$0=\beta\{(-2n+1)\lambda^2+(-\frac{3}{2}n-1)|v|^2\}$$
From our assumption we have $\beta\neq 0$, and consequently $\lambda=0$ and $v=0$. Substituting these into equation (\ref{eq:12}), we obtain $c=0$.
With the above-derived results in hand, we sum both sides of equation (\ref{eq:14}) over $i$ and obtain
$$2c(n-1)=\alpha\{-(n-1)(2\lambda^2+\frac{3}{2}|v|^2)\}+\beta\{\frac12|v|^2-\lambda^2-\sum_{i}|S_{ii}|^2\}$$
Substituting $\lambda=0$, $v=0$ and $c=0$ into the identity, we obtain $\beta\sum_{i}|S_{ii}|^2=0$. Since $\beta\neq 0$, it follows that $S_{ii}=0$. As $S$ is a diagonal matrix, we get $S=A+A^{\ast}=0$. By virtue of the lemma \ref{lemma3}, $g$ is K\"ahler flat.
\end{proof}
In addition to the case of the Riemannian connection, since the Bismut connection is an important connection in complex geometry, we also proposed a Bismut version of the problem, namely Theorem \ref{thm2}. We now give the proof of Theorem \ref{thm2} by an analogous method.
\begin{lemma}[\cite{ChenZ25}] \label{lemma4}
Given any Hermitian manifold $(M^n,g)$ and under any local unitary frame $e$, the components of Bismut curvature and Chern curvature are related by
\begin{equation}\label{eq:bismut-chern-rel}
\begin{cases}
T_{ik,\bar{j}}^{\ell}=R^{c}_{k\bar{j}i\bar{\ell}}-R^{c}_{i\bar{j}k\bar{\ell}},\\[4pt]
R^b_{i\bar{j}k\bar{\ell}}
= R^{c}_{i\bar{j}k\bar{\ell}}
+ T^\ell_{ik,\bar{j}}
+ \overline{T^k_{j\ell,\bar{i}}}
+ \sum_{s} \Bigl\{
T^s_{ik}\overline{T^s_{j\ell}}
- T^\ell_{is}\overline{T^k_{js}}
\Bigr\},\\[4pt]
\end{cases}
\end{equation}
where the subscripts after comma stand for covariant derivatives with respect to the Bismut connection $\nabla^b$.
\end{lemma}

Applying symmetrization to the second equation of Lemma \ref{lemma4}, we get
\begin{lemma}
\begin{equation}\label{eq:sym-hatRb-1}
4\widehat{R}^{b}_{i\overline{j}k\overline{l}}
=4\widehat{R}^{c}_{i\overline{j}k\overline{l}}
-\sum_{r=1}^{n}\Bigl(T^{j}_{ir}\overline{T}^{k}_{lr}
+T^{l}_{ir}\overline{T}^{k}_{jr}
+T^{j}_{kr}\overline{T}^{i}_{lr}
+T^{l}_{kr}\overline{T}^{i}_{jr}\Bigr)
\end{equation}
In particular, when \(i=j\neq k=l\), we have
\begin{equation}\label{eq:sym-hatRb-2}
4\widehat{R}^{b}_{i\overline{i}k\overline{k}}
=4\widehat{R}^{c}_{i\overline{i}k\overline{k}}
-\sum_{r=1}^{n}\Bigl(2\operatorname{Re}\bigl(T^{i}_{ir}\overline{T}^{k}_{kr}\bigr)
+\bigl|T^{k}_{ir}\bigr|^{2}+\bigl|T^{i}_{kr}\bigr|^{2}\Bigr)
\end{equation}
Meanwhile, we can also compute the following two identities.
\begin{equation}\label{eq:holo-sec-b}
R^{b}_{i\overline{i}i\overline{i}}
=R^{c}_{i\overline{i}i\overline{i}}
-\sum_{r=1}^{n}4\bigl|T^{i}_{ir}\bigr|^{2}
\end{equation}

\begin{equation}\label{eq:ric-b-def}
Ric_{i\overline{j}}^{b}
=\sum_{s=1}^{n}\bigl(R^{c}_{i\overline{s}s\overline{j}}
+R^{c}_{s\overline{j}i\overline{s}}
-R^{c}_{i\overline{j}s\overline{s}}\bigr)
\end{equation}
\end{lemma}

\begin{proof}[\textbf{Proof of Theorem \ref{thm2}}]
Let $(\mathfrak{g}, J, g)$ be a unimodular almost abelian Lie algebra with a Hermitian structure. Assume that the Bismut connection of $g$ has constant mixed curvature: $\mathcal{C}^{b}_{\alpha,\beta} = c$, where $\beta \neq 0$. Let  $e$ be an admissible frame and write $S=A+A^*$. After a unitary change of $\{e_2,...,e_n\}$ if necessary, we may assume that $S$ is a real diagonal matrix.
Using equations (\ref{eq:bismut-chern-rel}-\ref{eq:ric-b-def}) and (\ref{almostcur}), a direct computation gives
\begin{equation*}
\left\{
\begin{aligned}
R^{b}_{1\bar{1}1\bar{1}} &= -2\lambda^2 - \dfrac{3}{2}|v|^2, \quad R^{b}_{1\overline{1}i\overline{j}} = \big(-[A, A^*]+\lambda S\big)_{ij}, \quad R^{b}_{i\overline{j}1\overline{1}} = \big(-vv^*+S^2\big)_{ij},\\[4pt]
R^{b}_{i\overline{j}k\overline{l}} &= -S_{il}\overline{S}_{jk},\quad R^{b}_{i\overline{1}1\overline{l}} = \big([A, A^*]-\lambda S - S^2\big)_{il}, \quad \operatorname{Ric}^{b}_{1\overline{1}} = -3\lambda^{2}-2|v|^{2}, \quad \operatorname{Ric}^{b}_{i\overline{i}} = -|v_{i}|^{2}.
\end{aligned}
\right.
\end{equation*}
Here, the indices \(i, j, k, l\) are pairwise distinct, and \(i, j, k, l\) range from \(2\) to \(n\). Substituting the above identity into (\ref{eq:system-mixed}), we obtain the key identity below.
\begin{alignat}{2}
\label{eq:sys-line1} c&=\alpha\bigl(-3\lambda^{2}-2|v|^{2}\bigr)+\beta\bigl(-2\lambda^{2}-2|v|^{2}\bigr)\\[4pt]
\label{eq:sys-line2} c&=\alpha\bigl(-|v_{i}|^{2}\bigr)+\beta\bigl(-(S_{ii})^{2}\bigr)\\[4pt]
\label{eq:sys-line3} 2c&=\alpha\bigl(-3\lambda^{2}-2|v|^{2}-|v_{i}|^{2}\bigr)+\beta\bigl(-|v_{i}|^{2}+[AA^{*}]_{ii}-\lambda S_{ii}-(S^{2})_{ii}\bigr)
\end{alignat}
By adding equations \eqref{eq:sys-line1} and \eqref{eq:sys-line2} together, we get the following identity.
\begin{equation}\label{eq:combine-two-lines}
2c=\alpha\bigl(-3\lambda^{2}-2|v|^{2}-|v_{i}|^{2}\bigr)+\beta\bigl(-2\lambda^{2}-2|v|^{2}-(S_{ii})^{2}\bigr)
\end{equation}
Subtracting equation \eqref{eq:sys-line3} from equation \eqref{eq:combine-two-lines}, we obtain
\begin{align*}
0 &= \beta\left(2\lambda^{2}+2|v|^{2}+(S_{ii})^{2}-|v_{i}|^{2}+[A, A^{*}]_{ii}-\lambda S_{ii}-(S^{2})_{ii}\right)
\end{align*}
Since $\beta\neq 0$, we get
$$0=2\lambda^{2}+2|v|^{2}+(S_{ii})^{2}-|v_{i}|^{2}+[A, A^{*}]_{ii}-\lambda S_{ii}-(S^{2})_{ii},\quad \forall\,2\le i\le n.
$$
Summing over $i=2$ to $n$:
\begin{align*}
0=2(n-1)\lambda^{2}+2(n-1)|v|^{2}+\sum_{i=2}^{n}S_{ii}^{2}-|v|^{2}-\lambda\operatorname{tr}(S)-\operatorname{tr}(S^{2}).
\end{align*}
Thus
$$0=(2n-1)\lambda^{2}+(2n-3)|v|^{2}.$$
Since $n\ge 2$, we have $\lambda=0,\,v=0$.
By \eqref{eq:sys-line1}, $c=0$. Together with \eqref{eq:sys-line2}, this further implies $S_{ii}=0$ for all $i$.Recall that $S$ is diagonal and $S=A+A^*$, which implies $A+A^*=0$. Consequently, $g$ is K\"ahler flat.
\end{proof}

\section{Lie algebras with \texorpdfstring{$J$}{J}-invariant abelian ideal of codimension \texorpdfstring{$2$}{2}}
In this section, we consider a Lie algebra
${\mathfrak g}$ which contains a $J$-invariant abelian ideal ${\mathfrak a}$ of codimension $2$. Such a Lie algebra ${\mathfrak g}$ is generally not $2$-step solvable, but it is always solvable of step at most $3$. Throughout this section, we will follow the notations of \cite{GuoZ} and focus on the simpler case of $J{\mathfrak a} = {\mathfrak a}$.

Let $(\mathfrak{g},J,g)$ be a unimodular Lie algebra equipped with a Hermitian structure, and let ${\mathfrak a}\subset {\mathfrak g}$ be a $J$-invariant abelian ideal of codimension $2$. A unitary basis $\{ e_1, \ldots , e_n\}$ of ${\mathfrak g}^{1,0}$ will be called an {\em admissible frame} for $({\mathfrak g},J,g)$ if
$$ {\mathfrak a} = \mbox{span} \{ e_i+\overline{e}_i, \, \sqrt{-1}(e_i-\overline{e}_i); \ 2\leq i\leq n\},$$
As observed in \cite{GuoZ}, since ${\mathfrak a}$ is an abelian ideal, we have
\begin{equation*}
C^{\ast}_{ij}=D^j_{\ast i} =C^1_{\ast \ast} = D^{i}_{1\ast} =D^{\ast}_{1i} = 0, \ \ \ \ \ \forall \ 2\leq i,j\leq n,
\end{equation*}
and the only possibly non-zero components of $C$ and $D$ are
\begin{equation}  \label{CD-XYZ}
C^{j}_{1i}= X_{ij}, \ \ \ D^1_{11}=\lambda, \ \ \ D^j_{i1}=Y_{ij}, \ \ \ D^1_{ij}=Z_{ij}, \ \ \ D^1_{i1}=v_i, \ \ \ \ \ 2\leq i,j\leq n,
\end{equation}
where $\lambda \geq 0$, $v\in {\mathbb C}^{n-1}$ is a column vector, and $X$, $Y$, $Z$ are $(n-1)\times (n-1)$ complex matrices. Here and from now on, we assume that $D^1_{11}\geq 0$, which can be ensured by rotating $e_1$ by a complex number of norm $1$. Let $\varphi$ be the coframe dual to $e$, and for convenience write $\varphi=\,^t\!(\varphi_2, \ldots , \varphi_n)$. The structure equation then becomes
\begin{equation}
\left\{ \begin{aligned} &d\varphi_1 \, = \, -\lambda \varphi_1\overline{\varphi}_1 , \hspace{4.4cm} \\
&d\varphi \, = \, - \varphi_1\overline{\varphi}_1 \overline{v} -\varphi_1\,^t\!X \varphi + \overline{\varphi}_1 \overline{Y} \varphi - \varphi_1 \overline{Z} \,\overline{\varphi}. \end{aligned} \right.  \label{structure}
\end{equation}
Since $d^2\varphi =0$, we have
\begin{equation}
\left\{ \begin{aligned} \lambda (X^{\ast}\!+Y)+ [X^{\ast} ,Y]  -  Z\overline{Z} \, = \, 0, \\
\lambda Z - ( Z \,^t\!X + Y Z )\, =\, 0. \hspace{1.2cm} \end{aligned} \right.   \label{XYZ}
\end{equation}
Additionally, from (\ref{torsion}), we deduce that the only possibly non-zero components of the Chern torsion are as follows:
\begin{equation} \label{torsion-XYZ}
T^1_{1i}=v_i, \ \ \ \ T^1_{ij}=Z_{ji}-Z_{ij}, \ \ \ \ T^j_{1i}=Y_{ij}-X_{ij}, \ \ \ \ \ 2\leq i,j\leq n.
\end{equation}
From \cite[Proposition 3]{GuoZ}, we have
\begin{lemma} [\cite{GuoZ}]  \label{lemma6}
Let $({\mathfrak g},J,g)$ be a Lie algebra with Hermitian structure and let ${\mathfrak a}\subseteq {\mathfrak g}$ be a $J$-invariant abelian ideal of codimension $2$. Then under any admissible frame,
\begin{enumerate}
\item ${\mathfrak g}$ is unimodular $\Longleftrightarrow \, \lambda - \mbox{tr}(X)+\mbox{tr}(Y)=0$,
\item $g$ is K\"ahler  $\Longleftrightarrow v=0, \ ^t\!Z=Z, \text{and} \ X=Y$,
\item $g$ is Chern flat  $\Longleftrightarrow \lambda =0$, $v=0$, $Z=0$, $[Y,Y^{\ast}]=0$, and $[Y,X^{\ast}]=0$.
\end{enumerate}
\end{lemma}
Next, we recall the computational results of Li and Zheng in \cite{LZ25} for the expressions of Chern curvature components in the codimension two case. The following result follows from (16).
\begin{equation}\label{eq:16}
\left\{
\begin{aligned}
&\; R^{c}_{1\bar{1}1\bar{1}} = -2\lambda^2 - |v|^2,\quad R^{c}_{i\bar{i}i\bar{i}} = |Z_{ii}|^2,\quad \widehat{R}^{c}_{i\bar{i}k\bar{k}} = \dfrac{1}{4}|Z_{ik}+Z_{ki}|^2, \quad \forall\, 2\le i\ne k\le n,\\[6pt]
&\; \widehat{R}^{c}_{1\bar{1}i\bar{j}} = \dfrac{1}{4}\big(vv^* + [Y,Y^*] - \lambda(Y+Y^*) - Z\overline{Z} - {}^tZZ^* - {}^tZ\overline{Z}\big)_{ij}, \quad \forall\, 2\le i,j\le n,\\[6pt]
&\; \sum_{i=2}^n \widehat{R}^{c}_{1\bar{1}i\bar{i}} = \dfrac{1}{4}\big(|v|^2 - \lambda\operatorname{tr}(Y+Y^*) - 2\operatorname{tr}(Z\overline{Z}) - |Z|^2\big),
\end{aligned}
\right.
\end{equation}
Now we are ready to prove the bismut case of Theorem \ref{thm3} for $\mathfrak{g}$. Here again we do not need to assume that $\mathfrak{g}$ is unimodular.

Let $e$ be an admissible frame and write $B=Y-X$. A unitary change of basis $\{ e_2, \ldots , e_n\}$ via a unitary matrix $U$ transforms $Z$ into $UZ\,{}^t\!U$, and hence $\,^t\!Z+Z$ will be transformed into $U(^t\!Z+Z)\,^t\!U$. Since $\,^t\!Z+Z$ is a complex symmetric $(n-1)\times (n-1)$ matrix, it is well-known that there exists a unitary matrix $U$ such that $U(^t\!Z+Z)\,^t\!U$ is diagonal with real non-negative entries. Therefore, we may assume that $\,^t\!Z+Z$ is diagonal  with non-negative entries, and we will fix such an admissible frame from now on. Using
(\ref{eq:4}),(\ref{eq:5}) and (\ref{eq:16}), we obtain the following expressions:
\begin{equation}\label{cur-XYZ1}
\left\{
\begin{alignedat}{2}
& R_{1\bar{1}1\bar{1}} = -2\lambda^2 - \dfrac{3}{2}|v|^2, \quad R_{i\bar{i}i\bar{i}} = |Z_{ii}|^2 - \dfrac{1}{2}|B_{ii}|^2, \quad \forall\, 2\le i\le n, \\[6pt]
& \widehat{R}_{i\bar{i}k\bar{k}} = -\dfrac{1}{8}\big(|B_{ik}|^2 + |B_{ki}|^2 + 2\operatorname{Re}(B_{ii}\overline{B}_{kk})\big), \quad \forall\, 2\le i\ne k\le n, \\[6pt]
& \sum_{i=2}^{n} \widehat{R}_{1\bar{1}i\bar{i}} = \dfrac{1}{4}\big(|v|^2 - \lambda \operatorname{tr}(Y+Y^{*}) - 2\operatorname{tr}(Z\overline{Z}) - |Z|^2\big) - \dfrac{1}{8}\big(|v|^2 + \|B\|^2 + |{}^t Z - Z|^2\big).
\end{alignedat}
\right.
\end{equation}
where the last line is the result of taking trace in the previous line. We now turn our attention to the proof of Theorem \ref{thm1} for ${\mathfrak g}$.

\begin{proof}[{\bf Proof of Theorem \ref{thm3}}]
Let us fix a unimodular Hermitian Lie algebra  $(\mathfrak g, J, g)$ which contains an abelian ideal ${\mathfrak a}$ of codimension $2$ satisfying $J{\mathfrak a}={\mathfrak a}$. Assume that the Riemann connection of $g$ has constant mixed curvature: $\mathcal{C}_{\alpha,\beta} = c$, where $\beta \neq 0$.  Since $\mathcal{C}_{\alpha,\beta} = c$, by (\ref{eq:system-mixed}),(\ref{eq:lem-eq2}) and (\ref{cur-XYZ1}) We perform classification and obtain the following equation:
\begin{align}
\label{eq:syss-line1}
& c=\alpha\left\{-2\lambda^{2}-\frac{3}{2}|v|^{2}-\frac{1}{4}|Z+Z^{T}|^{2}\right\}
+\beta\left\{-2\lambda^{2}-\frac{3}{2}|v|^{2}\right\},\\[4pt]
\label{eq:syss-line2}
& c=\alpha\left\{-\frac12(BB^{*})_{ii}\right\}
+\beta\left\{|Z_{ii}|^{2}-\frac12|B_{ii}|^{2}\right\},\\[4pt]
\label{eq:syss-line3}
& 2c(n-1)
=\alpha\left\{
-(n-1)\left(2\lambda^{2}+\frac32|v|^{2}+\frac14|Z+Z^{T}|^{2}\right)
-\frac12|B|^{2}
\right\} \notag \\
&\quad +\beta\left\{
\frac12|v|^{2}+\lambda^{2}-\frac12|B|^{2}-|Z+Z^{T}|^{2}
\right\},\\[4pt]
\label{eq:syss-line4}
& 2cn(n-1)
=\alpha\bigl\{-(n-1)|B|^{2}\bigr\}+\beta\bigl\{|Z+Z^{T}|^{2}-\lambda^{2}-|B|^{2}\bigr\}.
\end{align}
Equation \eqref{eq:syss-line3} follows by summing the second identity in \eqref{eq:system-mixed} over $i=2$ to $n$, and Equation \eqref{eq:syss-line4} follows by summing the first identity in \eqref{eq:system-mixed} over $i\neq k$ with $i=2$ to $n$. When $\alpha=\beta$, since $\beta\neq 0$ from our hypothesis, dividing both sides of the above equation by $\beta$ yields the following simplified expression.
\begin{align}
\label{eq:syss-s-line1}
& P = -4\lambda^{2}-3|v|^{2}-\frac{1}{4}|Z+Z^{T}|^{2},\\[4pt]
\label{eq:syss-s-line2}
& 2P(n-1) = (-2n+3)\lambda^{2}+\left(-\frac{3}{2}n+2\right)|v|^{2}
+\left(-\frac{1}{4}n-\frac{3}{4}\right)|Z+Z^{T}|^{2}-|B|^{2},\\[4pt]
\label{eq:syss-s-line3}
& P = -\frac{1}{2}(BB^{*})_{ii}+|Z_{ii}|^{2}-\frac{1}{2}|B_{ii}|^{2},\\[4pt]
\label{eq:syss-s-line4}
& 2Pn(n-1) = -n|B|^{2}+|Z+Z^{T}|^{2}-\lambda^{2}.
\end{align}
where $P=\dfrac{c}{\beta}$. From equation \eqref{eq:syss-s-line4}, we can solve for \(|B|^2\):
\[
|B|^{2}=\frac{|Z+Z^{T}|^{2}-\lambda^{2}-2n(n-1)P}{n}.
\]
Substitute it into \eqref{eq:syss-s-line2}:
\[
2P(n-1)=(-2n+3)\lambda^{2}+(-\frac{3}{2}n+2)|v|^{2}
+(-\frac{1}{4}n-\frac{3}{4})|Z+Z^{T}|^{2}
-\frac{|Z+Z^{T}|^{2}-\lambda^{2}-2n(n-1)S}{n}.
\]
Multiply both sides by $n$:
\begin{align}
\label{eq:zero-identity}
0&=(-2n^{2}+3n+1)\lambda^{2}+\left(-\frac{3}{2}n^{2}+2n\right)|v|^{2}
+\left(-\frac{1}{4}n^{2}-\frac{3}{4}-1\right)|Z+Z^{T}|^{2}.
\end{align}
When $n\ge 2$,
\[
-2n^{2}+3n+1<0,\quad -\frac{3}{2}n^{2}+2n<0,\quad -\frac{1}{4}n^{2}-\frac{3}{4}-1<0.
\]
Hence equation \eqref{eq:zero-identity} yields
\[
\lambda^{2}=0,\quad |v|^{2}=0,\quad |Z+Z^{T}|^{2}=0
\quad\Longrightarrow\quad \lambda=0,\;v=0,\;Z+Z^{T}=0.
\]
From \eqref{eq:syss-line1}, we obtain $\dfrac{c}{\beta}=0$. Since $\beta\neq 0$, we get $c=0$.
Then $\mathrm{tr}(Z\overline{Z})=\lambda^{2}+\lambda\mathrm{tr}(Y+Y^{*})=0$. As $c=0$, equation \eqref{eq:syss-s-line2} yields $B=0$. The first identity in \eqref{XYZ} gives $[X^{*},Y]=0$. Since $B=Y-X$, we get $Y=X$. Hence $[X,Y]=0$ and $[X^{*},X]=0$. By Lemma \ref{lemma6}, we conclude that $g$ is K\"ahler flat.
\end{proof}

\begin{remark}
    \upshape  In the proof of Theorem 3 presented above, we have only established the conclusion under the additional assumption $\alpha=\beta$. Although we expect the assertion to hold for all $(\alpha,\beta)\neq (0,0)$, a complete proof for general parameters is still missing. We leave this general case for future investigation.
\end{remark}

\noindent\textbf{Acknowledgments.} We would like to express our sincere gratitude to Professor Fangyang Zheng for his insightful guidance, valuable discussions, and continued support throughout the course of this research.

\vs

\vspace{0.5cm}

\raggedbottom

\end{document}